\documentstyle [12pt,epsf,amssymb]{article}

\begin{document}
\def\l{\lambda}
\def\m{\mu}
\def\a{\alpha}
\def\b{\beta}
\def\g{\gamma}
\def\G{\Gamma}
\def\d{\delta}
\def\e{\epsilon}
\def\o{\omega}
\def\O{\Omega}
\def\v{\varphi}
\def\t{\theta}
\def\r{\rho}
\def\bs{$\blacksquare$}
\def\bp{\begin{proposition}}
\def\ep{\end{proposition}}
\def\bt{\begin{th}}
\def\et{\end{th}}
\def\be{\begin{equation}}
\def\ee{\end{equation}}
\def\bl{\begin{lemma}}
\def\el{\end{lemma}}
\def\bc{\begin{corollary}}
\def\ec{\end{corollary}}
\def\pr{\noindent{\bf Proof: }}
\def\note{\noindent{\bf Note. }}
\def\bd{\begin{definition}}
\def\ed{\end{definition}}
\def\C{{\mathbb C}}
\def\P{{\mathbb P}}
\def\Z{{\mathbb Z}}
\def\d{{\rm d}}
\def\deg{{\rm deg\,}}
\def\deg{{\rm deg\,}}
\def\arg{{\rm arg\,}}
\def\min{{\rm min\,}}
\def\max{{\rm max\,}}

\newtheorem{th}{Theorem}[section]
\newtheorem{lemma}{Lemma}[section]
\newtheorem{definition}{Definition}[section]
\newtheorem{corollary}{Corollary}[section]
\newtheorem{proposition}{Proposition}[section]

\centerline{\Large{\bf Analytic continuation and fixed points}}
\medskip
\centerline{\Large{\bf of the Poincar\'e mapping}}
\medskip
\centerline{\Large{\bf for a polynomial Abel equation}}
\par\bigskip
\bigskip
\bigskip
\centerline{\bf J.-P. Francoise} \centerline{Universit\'e de Paris
VI} \centerline{Laboratoire J.-L. Lions, UMR 7589 du CNRS}
\centerline{175 Rue de Chevaleret} \centerline{75013 Paris,
France}
\medskip
\centerline{\bf N. Roytvarf and Y. Yomdin} \centerline{Department
of Mathematics } \centerline{The Weizmann Institute of Science}
\centerline{Rehovot 76100, Israel}
\par\bigskip
\centerline{jpf@math.juissieu.fr}
\centerline{nina.roytvarf@weizmann.ac.il }
\centerline{yosef.yomdin@weizmann.ac.il }\bigskip
\bigskip
\vskip 10pt

The first author expresses his gratitude to the Weizmann Institute
for providing a financial support for a visiting during which this
article was initiated. The third author expresses his gratitude to
the University Paris VI and to the IHES for providing a financial
support for a visiting during which this article was continued.
The work of the second and the third authors has been supported
also by the ISF Grants No. 264-2002 and 979-2005, by the BSF Grant
No. 2002243 and by the Minerva Foundation.

\newpage
\centerline{\bf Abstract}
\bigskip
We consider an Abel differential equation $y'=p(x)y^2+q(x)y^3$
with $p(x)$, $q(x)$ -- polynomials in $x$. For two given points
$a$ and $b$ in $\mathbb C$, the ``Poincar\'e mapping" of the above
equation transforms the values of its solutions at $a$ into their
values at $b$. In this paper we study global analytic properties
of the Poincar\'e mapping, in particular, its analytic
continuation, its singularities and its fixed points (which
correspond to the ``periodic solutions" such that $y(a)=y(b)$). On
one side, we give a general description of singularities of the
Poincar\'e mapping, and of its analytic continuation. On the other
side, we study in detail the structure of the Poincar\'e mapping
for a local model near a simple fixed singularity, where an explicit solution can be
written. Yet, the global analytic structure (in particular, the
ramification) of the solutions and of the Poincar\'e mapping in
this case is fairly complicated, and, in our view, highly
instructive. For a given degree of the coefficients we produce
examples with an infinite number of {\it complex} ``periodic
solutions" and analyze their mutual position and branching. Let us
remind that Pugh's problem, which is closely related to the
classical Hilbert's 16th problem, asks for the existence of a
bound to the number of {\it real} isolated ``periodic solutions".
\vskip 3pt
\baselineskip 18pt
\par\bigskip
\newpage

\section {Introduction.}
\setcounter{equation}{0}

In this paper we start an investigation of the global analytic
properties of the ``Poincar\'e mapping'' $\phi$ for an Abel
differential equation of the form
$$ y' = p(x)y^2 + q(x)y^3. \ $$
For two given points $a$ and $b$, $\phi$ transforms the values
$y(a)$ of the solutions $y(x)$ of this equation at $a$ into their
values $y(b)$ at $b$. A more accurate definition is given in Section 3
below.

\medskip

We follow here the approach of the classical Analytic Theory of
Differential Equations (see \cite{Pain,Gol,Ince,Lor}) and consider
the Abel equation, its solutions and its Poincar\'e mapping in the
{\it complex domain}. Consequently, the main questions under
investigation are fixed and movable singularities, analytic
continuation and ramification of the functions involved.

\medskip

Below we always restrict ourselves to the case of a {\it
polynomial} Abel equation

\be y'=p(x)y^2+q(x)y^3,\quad y(0)=y_0 \ee with $p(x),q(x)$ --
polynomials in the complex variable $x$ of the degrees $d_1,d_2$
respectively.

The following two problems for Equation (1.1) have been
intensively studied (see \cite{Bfy1}-\cite{Bfy5}): \vskip 3pt

\noindent{\bf a.} For given $a,b$ is it possible to bound the
number of {\it real} solutions $y(x)$ of (1.1), satisfying
$y(a)=y(b)$, in terms of the degrees $d_1$, $d_2$ only? \vskip 3pt

\noindent{\bf b.} Is it possible to give explicit conditions on
$p$ and $q$ for $y(a)\equiv y(b)$ for all the solutions of (1.1)?
\vskip 3pt

\noindent These two problems are well known to be closely related
to the classical Hilbert's 16th problem (second part) and
Poincar\'e's Center-Focus problem for polynomial vector fields on
the plane.

Adopting the standard terminology in these problems, we shall call
solutions of (1.1), satisfying $y(a)=y(b)$, the ``closed" or the
``periodic" ones, and the conditions on $p,q, \ a,b$ for
$y(a)\equiv y(b)$ the center conditions.

Abel equations were first investigated and studied by Abel himself
as natural extensions of Ricatti equations. Abel found several
examples which are integrable (\cite{Ab}). Then this list was
enriched by Liouville. Classical references are
\cite{Ab,Liou,Gol,Ince}, and the modern references
\cite{Lor,Lli.Zhe.Wei,ChT.Roch} have been instrumental for us.
\vskip 3pt

The main motivation for our study comes from these classical
examples of the polynomial Abel equation which can be solved
explicitly. Moreover, we mostly (but not always) restrict
ourselves to one case where the first integral is rational.
However, the global analytic structure (in particular, the
ramification) of the solutions and of the Poincar\'e mapping in
this example turns out to be fairly complicated, and, in our
view, highly instructive. We study the singularities and the
branching of the solutions and of the Poincar\'e mapping. In
particular, for a fixed degree of the coefficients we produce,
varying the parameters, examples with arbitrarily many and with an
infinite number of ``periodic solutions". We analyze the mutual
position and branching of these periodic solutions.

\medskip

Our attempt to better understand the global structure of the
Poincar\'e mapping for the polynomial Abel equations was motivated
also by the recent progress in the investigation in
\cite{Bfy1}-\cite{Bfy5}, \cite{Bby,Bry,Pry} of the Center-Focus
problem for (1.1). As it was mentioned above, this problem is to
give explicit conditions on $p$ and $q$ for the Poincar\'e mapping
on $a,b$ to be identical, i.e. for $y(a)\equiv y(b)$. In
particular, in \cite{Bfy1}-\cite{Bfy5} the Moment and the
Composition conditions, providing a close approximation of the
center conditions, have been introduced. On this base in
\cite{Bby} the complete ``local center conditions'' have been
obtained, and the ``local Bautin ideal'' has been computed for the
Poincar\'e mapping $\phi$, while in \cite{Bry} similar conditions
``at infinity" have been found.

\medskip

Via Bautin's approach \cite{Ba}, further developed in
\cite{Fy1,Fy2,Y1} the knowledge of the Bautin ideal of the
Poincar\'e mapping $\phi$ allows one to produce ``semi-local''
bounds on the fixed points of $\phi$. In other words, we get a
fairly accurate control of the fixed points inside the disk of
convergence of the Taylor series of $\phi$ at the origin. Let us
remind that the problem of the {\it global control} of real fixed
points of $\phi$ is very closely related to the Hilbert 16th
problem of counting limit cycles of the plane vector-fields.

\medskip

However, the methods of \cite{Fy1,Fy2,Y1} are at present
absolutely limited to the disk of convergence of $\phi$. Any
attempt to ``globalize'' the information produced by these methods
will require a much better understanding of the global analytic
nature of $\phi$, in particular, of its analytic continuation, its
singularities and its ramification structure. In this paper we
start an investigation in this direction.

\medskip

The paper is organized as follows:

\medskip

In Section 2 we reprove the classical results of \cite{Pain} which
provide the description of singularities of the solutions of
(1.1). Our proof is somewhat more ``quantitative'' than the
classical one, providing an accurate estimates of the domains and
the parameters involved. We also prove some lemmas relating the
position of the singularities of the solutions of (1.1) with the
initial values of these solutions.
\par
On this base, in Section 3, we give an accurate definition of the
Poincar\'e mapping $\phi$, discuss the problem of the analytic
continuation of $\phi$, and give a constructive procedure of this
analytic continuation, based on the path deformation following the
moving singularities of the solutions. In Section 4 we describe
typical singularities of the Poincar\'e mapping $\phi$. This
completes our general description of the Poincar\'e mapping for
Abel equation.
\par
In Section 5 we discuss a local model of Abel equation near a simple fixed
singularity, and produce its explicit solutions.
\par
In Section 6 we analyze the singularities and the ramification of
the solutions and of the Poincar\'e mapping and we combine the results
to analyze the geometry of the periodic solutions.


\section{The Abel equation}
\setcounter{equation}{0}

Below we shall always assume that the functions $p(x)$ and $q(x)$
in the Abel equation (1.1) are polynomials in $x$ with complex
coefficients. Most of the results below remain valid for $p(x)$
and $q(x)$ - much more general analytic functions, but our
assumption simplifies a presentation.

\medskip

Let $a\in\mathbb C$. Denote by $y(y_a,x)$ the solution of the
equation (1.1), satisfying $y(y_a,a)=y_a$. By the uniqueness and
existence results for ordinary differential equations, the
solution $y(y_a,x)$ exists in a certain neighborhood of $a$ and is
there a regular complex analytic function of the complex argument
$x$. However, an analytic continuation of $y(y_a,x)$ may lead to
singularities.

\medskip

The classical result of Painlev\'e \cite{Pain} shows that the
``movable'' singularities of the solutions $y(y_a,x)$ must be
``algebroid''. Moreover, following the proof of Painlev\'e (see,
for example \cite{Pain,Gol,Lor}), one can easily show that at each
movable singular point $x_0,\ y(y_a,x)$ behaves as
${1\over\sqrt{x-x_0}}$. In order to relate singularities of $y$
with those of the Poincar\'e mapping $\phi$ we need more detailed
information on the position of singularities, on their dependence
on the initial values, etc., than is usually given. So we reprove
in the special case of the equation (1.1) the classical results,
providing all the required estimates.

\medskip

Notice that $y\equiv 0$ is a solution of (1.1). It follows, in
particular, that as $y_a\rightarrow 0$ all the singularities of
$y(y_a,x)$ tend to infinity. Below we make this remark more
precise.

\medskip

Another remark is that, as we shall see below, the problematic
points of the equation (1.1) are zeroes of $q(x)$. We denote these
zeroes $x_1,\cdots,x_m$ and always distinguish between the
``fixed" singularities of $y$ at $x_1,\cdots,x_m$ and the
``movable" singularities of $y$, which occur at points different
from $x_1,\cdots,x_m$.

\subsection{Domain of regularity of the solutions}

The following assumptions will be preserved along the rest of this
section: $p(x)$ and $q(x)$ are polynomials  of degree $m$ in $x$,
with $\Vert p\Vert,\ \Vert q\Vert\leq K$. The norm of a polynomial
is defined here as the sum of the absolute values of its
coefficients. Let $a\in \mathbb C$. Denote, as above, by
$y(y_a,x)$ the solution of the equation (1.1), satisfying
$y(y_a,a)=y_a$.

\bl Let $a\in {\mathbb C}, \ y_a\in \mathbb C$ be given. Then the
solution $y(y_a,x)$ exists in a disk $D_\rho(a)$ centered at $a$.
Here $\rho=\rho(\vert a \vert, \vert y_a\vert)$ is a positive
explicitly given function of its arguments, which for $\vert y_a
\vert$ big satisfies
$$ \rho(\vert a \vert, \vert y_a \vert)\geq C_1 (4K{\vert a
\vert}^m{\vert y_a \vert}^2)^{-1}.$$ For $\vert y_a \vert$ small
$\rho$ satisfies
$$ \rho(\vert a \vert, \vert y_a \vert)\geq
C_2 \left(\frac{1}{2K\vert y_a\vert}\right)^{\frac{1}{m+1}}.$$ In
particular, $\rho$ tends to infinity as $\vert y_a \vert$ tends to
zero.

\medskip

The solution $y(y_a,x)$ is bounded in the disk $D_\rho(a)$ by
$\hat y(\vert y_a\vert,\vert a\vert, \vert x-a\vert )$, with $\hat
y$ an explicitly given function of its arguments, satisfying
$$\hat y(\vert y_a\vert,\vert a\vert, t) \leq C_3(\vert y_a\vert,
\vert a\vert)(\rho - t)^{-{1\over 2}}.$$ \el

\pr For each $x\in {\mathbb C},\vert p(x)\vert\leq K{\vert x
\vert}^m,\vert q(x)\vert\leq K{\vert x \vert}^m$. Hence, the right
hand side of (1.1) is bounded in absolute value by $K{\vert x
\vert}^m(\vert y\vert^2+\vert y\vert^3)$. Therefore we get the
following differential inequality: \be \frac{d\vert y \vert}{dv}
\leq K{\vert x \vert}^m(\vert y\vert^2 +\vert y\vert^3). \ee  Here
$\frac{d}{dv}$ denotes a directional derivative at $x$ in any
(normalized) direction $v$ in the complex plane.

\medskip

Now consider a straight line $\ell$ in $\mathbb C$, passing
through $a$ and let $t$ be a normalized parameter along $\ell$,
with $t=0$ at $a$. Since for a running point $x(t) \in \ell, \
\vert x(t) \vert \leq \vert a\vert + t$ we get by (2.1) the
following differential inequality with respect to $t$: \be
\frac{d\vert y\vert}{dt} \leq K{(\vert a \vert + t)}^m(\vert
y\vert^2 +\vert y\vert^3). \ee Denote by $\tilde y(t)= \tilde
y(\vert a \vert, \vert y_a \vert,t)$ the solution of the
differential equation \be \frac{dy}{dt} = K{(\vert a \vert +
t)}^m(y^2 + y^3), \ee satisfying for $t=0$ the initial condition
$\tilde y(0)= \vert y_a \vert$. Then by (2.2) for each $t\geq 0$
we have $\vert y(y_a,x(t)) \vert \leq \tilde y(t)$.

\medskip

It remains to compute $\tilde y(t)$. Separating variables we
obtain $$\frac{dy}{y^2(y+1)}= dy({-1\over y}+{1\over
{y}^2}+{1\over {y+1}})=K{(\vert a \vert + t)}^m,$$ which gives
after integration the following implicit equation for the solution
$\tilde y(t)$, where we denote the function $ln(1+{1\over y}) -
{1\over y}$ by F(y): \be F(y) = F(\vert y_a \vert) + {K\over
{m+1}}({(\vert a \vert + t)}^{m+1}- {\vert a \vert}^{m+1}).\ee

The function $F(y)$ for $y$ positive is a negative strictly
increasing function, tending to $-\infty$ as $y$ tends to zero,
and approaching zero from below as $y$ tends to $\infty$. The
inequality $\frac{1}{y^2(y+1)}< \frac{1}{y^3}$ shows that $F(y)>
-\frac{1}{2y^2}$. To bound $F(x)$ from above, define
$h(y)=\frac{1}{2y^3}$ for $y\geq 1$, and $h(y)=\frac{1}{2y^2}$ for
$y\leq 1$. We have $h(y)< \frac{1}{y^2(y+1)}$, and therefore
$F(x)<H(y)$, where $H(y)= \int h(y)dy = -\frac{1}{4y^2}$ for
$y\geq 1$ and $H(y)= -\frac{1}{2y}+\frac{1}{4}$ for $y \leq 1$.
Finally we get the following bounds from two sides for $F(y)$: \be
-\frac{1}{2y^2} < F(y) < H(y).\ee Solving the equation (2.4) and
taking into account the lower bound in (2.5) we get $$ \tilde y(t)
< [-2(F(\vert y_a \vert) - {K\over {m+1}}({(\vert a \vert +
t)}^{m+1}- {\vert a \vert}^{m+1}))]^{-{1\over 2}}.$$ Now applying
the upper bound in (2.5) we finally obtain the following
inequality: \be \tilde y(t) < \hat y(t) = [-2(H(\vert y_a \vert) -
{K\over {m+1}}({(\vert a \vert + t)}^{m+1}- {\vert a
\vert}^{m+1}))]^{-{1\over 2}}.\ee Taking into account an explicit
definition of the function $H$ given above, we obtain $$ \hat y(t)
= [\frac{1}{4{\vert y_a \vert}^2} - {K\over {m+1}}({(\vert a \vert
+ t)}^{m+1}- {\vert a \vert}^{m+1}))]^{-{1\over 2}}, \ \vert y_a
\vert \geq 1,$$ and \be \hat y(t) = [\frac{1}{2{\vert y_a \vert}}
- \frac{1}{4} - {K\over {m+1}}({(\vert a \vert + t)}^{m+1}- {\vert
a \vert}^{m+1}))]^{-{1\over 2}}, \ \vert y_a \vert \leq 1.\ee The
function $\hat y(t) = \hat y(\vert a \vert, \vert y_a \vert,t)$
grows with $t$. For $t=0$ it takes value $2\vert y_a \vert$, if
$\vert y_a \vert \geq 1$, and it takes value $({\frac{1}{2{\vert
y_a \vert}} - \frac{1}{4}})^{-{1\over 2}}$, if $\vert y_a \vert
\leq 1$ (which is of order $\sqrt {2\vert y_a \vert}$ for small
$\vert y_a \vert$). This function is finite until the expression
in the parentheses remains positive. This gives the following
expression for the radius of the disk of existence of the
solution: \be \rho(\vert a \vert, \vert y_a \vert)= \vert a
\vert([1-{{m+1}\over {K{\vert a \vert}^{m+1}}}H(\vert y_a
\vert)]^{1\over {m+1}}-1). \ee For $\vert y_a \vert$ big this
gives us the following expression: \be \rho(\vert a \vert, \vert
y_a \vert)\simeq (4K{\vert a \vert}^m{\vert y_a \vert}^2)^{-1}.\ee
For $\vert y_a \vert$ small we get \be \rho(\vert a \vert, \vert
y_a \vert)\simeq \left(\frac{1}{2K\vert
y_a\vert}\right)^{\frac{1}{m+1}}. \ee This completes the proof of
Lemma 2.1.

\bc $y(y_a,x)$ is regular in the disk $D_R$ of radius $R$ centered
at the origin, with $R$ growing as $\left(\frac{1}{2K\vert
y_a\vert}\right)^{\frac{1}{m+1}}$, as $\vert y_a\vert$ tends to
zero.\ec

\subsection{Singularities of $y(y_a,x)$}

In this subsection we reprove the classical results on the
structure of singularities of $y(y_a,x)$, stressing the explicit
estimates of the size of the domains, where the results are valid.
The assumptions on $p$ and $q$ and the notations remain the same
as in subsection 2.1.

\bl If $x_0\in\mathbb C$ is a singular point of the solution
$y(x)$ of (1.1), then $y$ tends to infinity as $x$ tends to $x_0$.
\el \pr If $\lim_{x\rightarrow x_0}\ y(x)\not = \infty$, then
there is a constant $c>0$ and a sequence $x_i$ converging to $x_0$
such that $\vert y(x_i) \vert\leq c$ for each $i$. Applying Lemma
2.1 to one of the points $x_i$, taken as $a$, we obtain that
$y(x)$ is regular inside a disk around $x_i$, of a certain radius
$\rho > 0$ which does not depend on $i$. Taking $x_i$ sufficiently
close to $x_0$, we conclude that $y$ is regular at and around
$x_0$. This contradiction proves the lemma.

\medskip

Now we give an analytic description of the movable singular points
of the solutions of (1.1). As we have already mentioned in the
Introduction, the classical result of Painlev\'e
\cite{Pain,Gol,Lor} shows that the ``movable'' singularities of
the solutions $y(y_a,x)$ of the equation (1.1) must be
``algebroid''. Moreover, following the proof of Painlev\'e (see,
for example \cite{Pain,Gol,Lor}), one can easily show that at each
movable singular point $x_0,\ y(y_a,x)$ behaves as
${1\over\sqrt{x-x_0}}$. However, in order to relate singularities
of $y$ with those of the Poincar\'e mapping $\phi$ we need more
accurate estimates, than are usually given, and in particular, we
have to describe the behavior of the movable singular points of
the solutions of (1.1) as the function of the initial value $y_a$.
So we reprove below the classical result of Painlev\'e (in the
special case of the equation (1.1)), providing the required
estimates.

\medskip

To simplify the statement of the results below, let us introduce
some notations. Assume that $x_0\in \mathbb C$ is given, $x_0$
different from the zeroes $x_1,\cdots,x_m$ of $q(x)$. Let
$\eta(x_0)=\vert q(x_0)\vert>0$. Put $R(x_0)=2(\vert x_0\vert+1)$
and define $r(x_0)=min\left({1\over 4}R,\frac{\eta}{2m(K+1)
R^{m-1}}\right)$, where $K$ as above is the maximum of the norms
$\Vert p\Vert,\Vert q\Vert$. Finally, let us define
$M(x_0)=\frac{4m(K+1)R^m}{\eta}$ and put $\delta(x_0)= {1\over M}$
and $c(x_0)={4\over {M\eta}}=\frac{1}{m(K+1) R^m}$.

\bt For any $x_0\in {\mathbb C}, \ x_0$ different from the zeroes
$x_1,\cdots,x_m$ of $q(x)$, there is a unique solution
$y(x)=y^{[x_0]}(x)$ of (1.1) with a singularity at $x_0$. This
solution has an algebraic ramification of order 2 at $x_0$. In a
neighborhood of its singular point $x_0$, the solution
$y^{[x_0]}(x)$ of (1.1) is given by the Puiseux series

\be y^{[x_0]}(x)=\frac{c(x_0)}{(x-x_0)^{1/2}}
\cdot(1+\sum^\infty_{k=1} c_k(x_0)(x-x_0)^ {\frac{k}{2}}),\ee
converging for $\vert x-x_0\vert \leq r(x_0)$, with the
coefficients $c(x_0), \ c_k(x_0)$ univalued analytic functions in
${\mathbb C} \setminus \{x_1,\cdots,x_m \}$, satisfying there
$\vert c_k(x_0) \vert \leq \delta (x_0) {(r(x_0))}^{-k}.$ \et

\pr The proof of Theorem 2.1 takes the rest of the present
subsection. In the course of the proof we stress some intermediate
steps which will be used later as lemmas, propositions, etc.
\par
After the change of the dependent variable $u=\frac{1}{y}$ the
equation (1.1) takes the form \be \frac{du}{dx} =
-p(x)-q(x)\frac{1}{u} = -\frac{p(x)u+q(x)}{u}. \ee Changing the
roles of the dependent and the independent variables, we get \be
\frac{dx}{du} = - \frac{u}{p(x)u+q(x)} = G(u,x). \ee The right
hand side of (2.13) is a regular function of $u$ and $x$ near
$(0,x_0 )$, since $q(x_0)\not = 0$. We need an upper bound on
$G(u,x)$ in an explicitly given neighborhood of $(0,x_0 )$. Let us
remind that in order to simplify the statement of this and further
bounds, we have introduced the following notations: $\vert
q(x_0)\vert=\eta>0$, $R=2(\vert x_0\vert+1)$, $r=min\left({1\over
4}R,\frac{\eta}{2m(K+1) R^{m-1}}\right)$, where $K$ as above is
the maximum of the norms $\Vert p\Vert,\Vert q\Vert$. Finally,
$M=\frac{4m(K+1)R^m}{\eta}$, $\delta = {1\over M}$ and $c={4\over
{M\eta}}=\frac{1}{m(K+1) R^m}$.

\bp For $\vert u\vert\leq\delta$ and $\vert x-x_0\vert\leq r$ we
have $\vert G(u,x)\vert \leq c$.\ep

\pr Inside the disk $D_R$ centered at the origin the derivative
$q'(x)$ of the polynomial $q(x)$ satisfies $\vert q'(x)\vert\leq
mKR^{m-1}$. By the choice of $r$ we have $D_r(x_0)\subset D_R$.
Hence for $x\in D_r (x_0)$,
$$ \vert q(x)\vert\geq \vert q(x_0)\vert
-r\cdot mKR^{m-1}\geq\frac{1}{2}\eta, $$ by the choice of $r$. On
the other hand, inside $D_R$ the inequality $\vert p(x)\vert\leq
KR^m$ is satisfied, and for $\vert u\vert\leq
\delta=\frac{1}{M}=\frac{\eta}{4m(K+1)R^m}$, we have $\vert p(x)u
\vert\leq\frac{1}{4}\eta$. Hence the absolute value of the
denominator $p(x)u +q(x)$ of $G(u,x)$ is at least
$\frac{1}{4}\eta$ for $x\in D_r(x_0)$ and $\vert
u\vert\leq\delta$. Under the same assumptions on $x$ and $u$ we
finally obtain
$$ \vert G(u,x)\vert\leq\frac{\vert u\vert}{\vert p(x)u+q(x)\vert} \leq\frac{
1/M}{(1/4)\eta} = \frac{1}{m(K+1) R^m}=c.$$ This completes the
proof of the proposition.

\medskip

Returning to the proof of Theorem 2.1 we see that for $\vert
u\vert\leq \delta$ and $x\in D_r(x_0)$ the following differential
inequality is satisfied:

\be \frac{d\vert x(u)\vert}{dv} \leq c \ee in any normalized
direction $v$. Hence for any $u$ with $\vert u\vert\leq\delta$
(and assuming that $x$ remains in $D_r(x_0)$) we obtain
$$ \vert x(u)-x_0\vert\leq c\delta =\frac{1}
{m(K+1)R^m}\cdot\frac{\eta}{4m(K+1)R^m} = $$
$$ = \frac{\eta}{4(m(K+1)R^m)^2} < \frac{\eta}{4m(K+1)R^{m-1}}=\frac{1}{2}r.$$
Hence for $\vert u\vert\leq\delta$ the solution $x(x_0,u)$ of
(2.13), satisfying $x(0) =x_0$, exists and it indeed remains in
the disk $D_r(x_0)$. This justifies {\it a posteriori} the use of
the differential inequality (2.14).

\medskip

Therefore $x(x_0,u)$ is a regular analytic function of two complex
variables $u$ and $x_0$ defined in the domain $\O=\{x_0 \in
{\mathbb C} \setminus \{x_1,\cdots,x_m \}, \ \ \vert u\vert \leq
\delta(x_0)\}.$ Consequently, $x(x_0,u)$ can be represented by a
converging power series \be
x(x_0,u)=A(x_0)+u\sum^\infty_{k=0}A_k(x_0)u^k, \ee with $A,A_k$
analytic functions of $x_0$, univalued and regular in ${\mathbb C}
\setminus \{x_1,\cdots,x_m \}.$

\medskip

Let us show that in fact

\be x(x_0,u)=x_0+u^2\sum^\infty_{k=0}a_k(x_0)u^k, \ \ a_0(x_0)\neq
0 \ee Indeed, the initial condition $x(x_0,u)=x_0$ implies
$A(x_0)\equiv x_0$. Next, the equation (2.13) shows that the
derivative $\frac {dx}{du}$ vanishes for $u=0$ and hence there are
no linear in $u$ terms in (2.15). A direct computation shows that
$a_0(x_0)=-\frac{1}{2q(x_0)}$ and hence $\vert a_0(x_0)\vert =
{1\over {2\eta(x_0)}}.$ On the other hand, as it was shown above,
for $\vert u\vert\leq\delta(x_0)$ the solution $x(x_0,u)$ of
(2.13), satisfying $x(0) =x_0$, exists and it remains in the disk
$D_r(x_0)$. Hence by the Cauchy formula for the Taylor
coefficients we get $\vert a_k(x_0)\vert \leq
r(x_0){\delta(x_0)}^{-(k+2)}.$

\medskip

Now the standard manipulations with the power series show that the
solution $u(x)=u^{[x_0]}(x)$ of the equation (2.12), satisfying
$u^{[x_0]}(x_0)=0$, as a function of $x$ has a ramification of
order 2 at $x_0$, and it can be represented by a Puiseux series

\be u(x)=u^{[x_0]}(x)=\sum^\infty_{k=1}
b_k(x_0)(x-x_0)^{\frac{k}{2}},\ b_1(x_0)\not = 0,\ee converging in
the disk $\vert x-x_0\vert\leq\xi(x_0)$, where $\xi(x_0)$ can be
given explicitly through the parameters defined above. This
completes the proof of Theorem 2.1. Note that the coefficients
$b_k(x_0)$ are univalued and regular functions of $x_0$ in
${\mathbb C} \setminus \{x_1,\cdots,x_m \}.$

\bc In a neighborhood of its singular point $x_0$, the solution
$y(x)=y^{[x_0]}(x)$ of (1.1) is given by

\be y(x)=y^{[x_0]}(x)=\frac{c(x_0)}{(x-x_0)^{1/2}}
\cdot(1+\sum^\infty_{k=1} c_k(x_0)(x-x_0)^ {\frac{k}{2}}),\ee
converging for $\vert x-x_0\vert \leq\xi_1(x_0)$, with
$c(x_0)={1\over b(x_0)}$.\ec

\pr $$y(x)=\frac{1}{u(x)}=$$ $$ = \frac{1}{b(x_0)(x-x_0)^{1/2}}
\cdot{1\over{1+\sum^\infty_{k=1}{{b_{k+1}(x_0)}\over
{b_1(x_0)}}(x-x_0)^{k/2}}}=$$
$$ = \frac{1}{b(x_0)(x-x_0)^{1/2}} \left( 1+\sum^\infty_{k=1} c_k(x_0)
(x-x_0)^{k/2}\right) . $$

\subsection{Singularities of $y(y_a,x)$ as functions of $y_a$}

In order to relate the singularities of the Poincar\'e mapping
with those of the solutions of (1.1), it is important to see how
the initial value of the solution $y$ at a certain regular point
influences the position of the singularities of $y$. The
description of the singularities of $y$ given above, allows one to
get a rather accurate information in this respect.
\par
Let us fix a certain point $c\in {\mathbb C},\ c\not=
x_1,\cdots,x_m$ (i.e. $q(c)\not = 0$).

\bl For any $y_c$, sufficiently large in absolute value, the
solution $y(y_c,x)$ of (1.1), satisfying $y(y_c,c)=y_c$, has a
singularity $x_0=x_0(y_c)$ in a neighborhood of $c$. The position
of this singularity, $x_0(y_c)$, is a regular function of $y_c$
for $\vert y_c\vert$ sufficiently large, and
$\frac{dx_0}{dy_c}\not=0$.\el

\pr Rewrite the expression (2.16) representing the solution
$x(x_0,u)$ of the equation (2.13), satisfying $x(x_0,0)=x_0$, as
$x(x_0,u)=x_0+u^2 g(x_0,u).$ Here $g(x_0,u)$ is a regular function
in in the domain $\O=\{x_0 \in {\mathbb C} \setminus
\{x_1,\cdots,x_m \}, \ \ \vert u\vert \leq \delta(x_0)\}$ defined
above, and $g(x_0,0)\neq 0.$ Now let us require that $x(x_0,u)=c$,
with $c$ fixed. We get an equation \be c=x_0+u^2 g(x_0,u) =
G(x_0,u) \ee between the position $x_0$ of the zero of $u$ (i.e.
of the singularity of $y={1\over u})$ and the value of $u$ at $c$,
$u=u(c)={1\over y(c)}$. Now for $u=0$ and $x_0=c,\ \frac{\partial
G}{\partial x_0}=1$. Hence by the implicit function theorem, $x_0$
is a regular function of $u=u(c)$, with $x_0(0)=c$. Moreover,
$\frac{dx_0}{du} = -{{\frac{\partial G}{\partial u}}\over
{\frac{\partial G}{\partial x_0}}} \neq 0$ for small $u\not= 0$,
since $\frac{\partial G}{\partial u}\not= 0$ for such $u$. This
completes the proof.

\medskip

\noindent {\bf Remark.} \ \ Explicit constants can be given in the
statement of Lemma 2.3, in the same terms as above.

\bc Let the solution $y(y_a,x)$ of (1.1), satisfying
$y(y_a,a)=y_a$ and continued to $x_0\not= x_1,\cdots,x_m$ along a
path $s$ in ${\mathbb C}\setminus \{x_1,\cdots,x_m\}$, have a
singularity at $x_0$. Then $x_0$ is a regular function of the
initial value $y_a$, and $\frac{dx_0}{dy_a}\not= 0$.\ec

\pr $y(y_a,x)$ tends to $\infty$ as $x$ tends to $x_0$. Fix a
point $c$ on the path $s$, sufficiently close to $x_0$, such that
for $y_c= y(y_a,c)$ and $c$ the conditions of lemma 2.3 are
satisfied. By this lemma, the singular point $x_0$ of the solution
is a regular function of the initial value $y_c$, and
$\frac{dx_0}{dy_c}\not= 0$. We have $\frac{ dy_c}{dy_a}\not=0$,
the mapping $y_a\rightarrow y_c$ being the Poincar\'e mapping of
(1.1) along the path $s$. We obtain that $x_0$ is a regular
function of $y_a$, with $\frac{dx_0}{dy_a} =
\frac{dx_0}{dy_c}\cdot \frac{dy_c}{dy_a}\not= 0$.

\medskip

Using the equation (2.19) we can extend Lemma 2.3 above and
describe the dependence of the position $x_0$ of the singular
point of the solution on the value of this solution at the
``original singular point" itself. Let us fix, as above, a certain
point $c\in {\mathbb C},\ c\not= x_1,\cdots,x_m$ (i.e. $q(c)\not =
0$). \bp For $y_c$ near $\infty$,  the position $x_0(y_c)$ of the
singularity of the solution $y(y_c,x)$ of (1.1), satisfying
$y(y_c,c)=y_c$, can be represented by a convergent Taylor series
in $u=u(c)={1\over {y(c)}}$ \be x_0-c = u^2\sum_{k=0}^{\infty}
\alpha_k u^k, \ \ \alpha_0 \neq 0. \ee Conversely, the value
$y(y_c,c)=y_c$ at $c$ of the solution $y$ of (1.1) having
singularity at $x_0$, can be represented by a convergent
fractional Puiseux series \be u=u(c)=\sum_{k=1}^{\infty} \beta_k
(x_0-c)^{k\over 2}, \ \ \beta_1 \neq 0.\ee \ep \pr We rewrite the
equation (2.19) in the form \be x_0-c=-u^2 g(x_0,u) = -G(x_0,u)
\ee between the position $x_0$ of the zero of $u$ (i.e. of the
singularity of $y={1\over u})$ and the value of $u$ at $c$,
$u=u(c)={1\over y(c)}$. Now, as above, for $u=0$ and $x_0=c,\
\frac{\partial G}{\partial x_0}=1$. Hence by the implicit function
theorem, $x_0-c$ is a regular function of $u=u(c)$, with
$x_0(0)-c=0$. Therefore, $x_0-c$ can be represented by the
convergent Taylor series \be x_0-c = u^2\sum_{k=0}^{\infty}
\alpha_k u^k.\ee Solving the equation (2.23) with respect to $u$
we get \be u=u(c)=\sum_{k=1}^{\infty} \beta_k (x_0-c)^{k\over
2}.\ee This completes the proof of the proposition. Below we shall
use it to describe the ramification of the Poincar\'e mapping
around its singularities.

\section{Analytic continuation of the Poincar\'e mapping}
\setcounter{equation}{0}

In this section we give an accurate definition of the Poincar\'e
mapping $\phi$ of the Abel Equation
$$ y\prime = p(x)y^2 + q(x)y^3 \ \leqno(1.1) $$
and discuss some problems related to the investigation of $\phi$.
Then we give a ``semi-constructive'' description of the analytic
continuation of $\phi$ along a given path.
\par
Let $a,b\in {\mathbb C},\ b\not= x_1,\cdots,x_m,$ where
$x_1,\cdots,x_m$ are, as above, all the zeroes of $q$. Notice that
if a solution $y$ of (1.1) happens to have a singularity at one of
the $x_i, \ i=1,\dots,m$, the analytic structure of this solution
near $x_i$ may be much more complicated than that described in
Section 2 above (see examples in Section 5 below). This is because
the equation (2.13) has a singularity at $(0,x_i)$; both the
numerator and the denominator of $G(u,x)$ vanish at this point.

\medskip

Let $s$ be a path in $\mathbb C$, joining $a$ and $b$. We {\it do
not assume} that $s$ avoids the points $x_1,\cdots,x_m$, unless
specifically stated. Let the initial value $y^0_a \in {\mathbb C}$
be given. Assume that the solution $y(y_a^0,x)$ of the equation
(1.1) satisfying $y(y^0_a,a)=y^0_a$ can be analytically continued
along $s$ from a neighborhood of $a$. In particular, this
continuation does not have singularities on $s$. \bd The (germ at
$y^0_a$ of the) Poincar\'e mapping $\phi=\phi_{a,b,s,y^0_a}$ of
the equation (1.1) along the path $s$ is defined as follows: it
associates to each $y_a$ near $y^0_a$ the value $y_b$ at the point
$b$ of the solution $y(y_a,x)$ of (1.1), satisfying
$y(y_a,a)=y_a$, and continued to $b$ along $s$. Thus
$\phi(y_a)=y(y_a,b)=y_b$.\ed Since $y\equiv 0$ is a solution of
(1.1), the germ of the Poincar\'e mapping at zero satisfies
$\phi(0)=0$ along any path $s$ and for any endpoints $a,b$.
Moreover, by Corollary 2.1, for any $R>0$ the solutions $y(y_a,x)$
are regular in the disk $D_R$, assuming that $\vert y_a\vert$ is
sufficiently small. Hence, for any $a,b,s$, the germ at the origin
of $\phi_{a,b,s}$ is defined and it does not depend on the path
$s$.

However, for larger values of $y_a$, the analytic continuation of
$y(y_a,x)$ along different paths $s$ may lead to different values
of $y_b$.

\medskip

Now assume that a path $\sigma$ from $w_0=\sigma(0)$ to
$w_1=\sigma(1)$ in the plane of the initial values $y_a$ is given,
parametrized by [0,1]. Assuming that none of the solutions
$y(w_t,x), \ w_t=\sigma(t), \ t \in [0,1],$ has a singularity on
the path $s$, the definition above works well and defines the
values (in fact, the germs) of $\phi(w_t),\ t\in[0,1]$, and, in
particular, $\phi(w_1)$,
\par
The problem appears if the singularities of the solutions
$y(w_t,x)$, continued along $s$, approach and cross the path $s$.
The idea of the following construction is that if we can deform
the path $s$ (following the movement of $w_t$ along the curve
$\sigma$) in such a way that it escapes the singularities of
$y(w_t,x)$, we can still  use this deformed path for the analytic
continuation of the solutions $y(w,x)$, and hence for the analytic
continuation of $\phi$.

\medskip

Let $\sigma$ as above be given. Assume that there exists a family
$s_t,\ t\in[0,1]$, of the paths from $a$ to $b$, with the
following properties:
\begin{itemize}
\item[1.]$s_t$ is a continuous in $t$ deformation of the original
path $s,\ s_0=s$.
\item[2.]For each $t\in[0,1]$, the solution $y(w_t,x)$, continued
along $s_t$, is regular at each point of $s_t$.
\end{itemize}

\bt The germ of the Poincar\'e mapping along $s$ at the point
$w_0$, $\phi_{a,b,s,w_0}$, allows for the analytic continuation
along $\sigma$ from $w_0=\sigma(0)$ to $w_1=\sigma(1)$. The
continued germ at $w_1$ is $\phi_{a,b,s_1,w_1}$ provided by the
analytic continuation of the solutions starting near $w_1$ along
$s_1$.\et \pr We shall show that for each $t\in [0,1]$ the value
(the germ) of $\phi(w_t)$, obtained by an analytic continuation of
$\phi$ along $\sigma$, is given by \be \phi(w_t)=y(w_t,b), \ee
with the solution $y(w_t,x)$, satisfying $y(w_t,a)=w_t$, being
analytically continued from $a$ to $b$ along the path $s_t$. We
can subdivide the process of the analytic continuation of $\phi$
into a finite number of small successive steps. In each step we
first move $w_t$ along $\sigma$, {\it not deforming} $s_t$
(providing that the singularities of $y(w_t,x)$ do not hit $s_t$).
Then we deform $s_t$, not changing $w_t$. Clearly, the first part
of each step gives an analytic continuation of $\phi$ along
$\sigma$, while the second step does not change $\phi$ at all.
Therefore the total procedure provides the required analytic
continuation of $\phi$. This completes the proof of Theorem 3.1.

\medskip

Theorem 3.1 reduces the problem of the analytic continuation of
the Poincar\'e mapping $\phi$ to a construction of the family of
paths $s_t$, with the properties stated above.

One particular case is very easy: if the singularities of
$y(w_t,x)$ do not approach the path $s$, it does not need to be
deformed, and we can take $s_t\equiv s$.

\medskip

\section{Singularities of the Poincar\'e mapping}
\setcounter{equation}{0}

Let $s$ be a path in $\mathbb C$ joining two points $a$ and $b$,
and let the initial value $y^0_a \in {\mathbb C}$ be given. Assume
that the solution $y(y_a^0,x)$ of the equation (1.1) satisfying
$y(y^0_a,a)=y^0_a$ can be analytically continued along $s$ from a
neighborhood of $a$ to each point of $s$ except, possibly, the
endpoint $b$. In particular, this continuation does not have
singularities in the interior points of $s$. If $b$ is also a
regular point of this solution, then the germ at $y^0_a$ of the
Poincar\'e mapping $\phi=\phi_{a,b,s,y^0_a}$ along the path $s$ is
defined and regular.

Consider now the case when the analytic continuation along $s$ of
the solution $y(y_a^0,x)$ has a singularity at $b$. From now on, we
assume that $b$ is different from the fixed singularities
$x_1,\dots,x_m.$ \bp Under the above assumptions there is a germ
of a real curve $\gamma \subset {\mathbb C}$ at $y^0_a$ such that
for $y_a \not \in \gamma$ the analytic continuation along $s$ of
the solution $y(y_a,x)$ is regular at each point of $s$ including
the endpoint $b$ \ep \pr We are in a situation of Corollary 2.3
above. By this corollary, the position $x_0(y_a)$ of the
singularity of the solution $y(y_a,x)$ near $b$ is a regular
function of $y_a$ near $y_a^0$. So the curve $\gamma$ is formed by
exactly those $y_a$ for which this singularity $x_0(y_a)$ belongs
to $s$. By the description of the singularities of the solutions
of (1.1) given in Section 2, $x_0(y_a)$ is the only singular point
of the local branch of the solution $y(y_a,x)$ near $b$. On the
other side, by the assumptions, there are no singularities of the
solution $y(y_a,x)$ on $s$ not in a neighborhood of $b$.
Therefore, for $y_a \not \in \gamma$ the analytic continuation
along $s$ of the solution $y(y_a,x)$ is regular at each point of
$s$, including the endpoint $b$. This completes the proof of the
proposition.

Now we are ready to describe the generic singular points of the
Poincar\'e mapping. \bt Let $s$ be a path in $\mathbb C$ joining
two points $a$ and $b \neq x_1,\dots,x_m$, and let the initial
value $y^0_a \in {\mathbb C}$ be given. Assume that the solution
$y(y_a^0,x)$ of the equation (1.1) satisfying $y(y^0_a,a)=y^0_a$
can be analytically continued along $s$ from a neighborhood of $a$
to each point of $s$ except the endpoint $b$, where this solution
has a singularity. Then for each $y_a$ in a neighborhood of
$y^0_a$, such that $y_a \not \in \gamma$, where the curve $\gamma$
has been defined in Proposition 4.1, the germ at $y_a$ of the
Poincar\'e mapping $\phi=\phi_{a,b,s,y_a}$ along the path $s$ is
defined and regular. In a punctured neighborhood $U_0$ of $y^0_a$
these germs can be analytically continued across $\gamma$, to form
a double-valued regular in $U_0$ function $\phi$, which allows for
a representation by a  Puiseux series \be y(b) = {1 \over {\sqrt
{(y_a-y_a^0)}}} \sum_{k=0}^{\infty} \nu_k (y_a-y_a^0)^{k\over 2},
\ \ \nu_0 \neq 0, \ee convergent in $U_0$. \et \pr The fact that
for each $y_a$ in a neighborhood of $y^0_a$, such that $y_a \not
\in \gamma$, the germ at $y_a$ of $\phi$ is defined and regular,
follows directly from Proposition 4.1. A possibility of the
analytic continuation of $\phi$ across $\gamma$ follows from the
results of Section 3. The local form of $\phi$ near $y^0_a$ can be
obtained in two ways.

The first one uses Corollary 2.1 and Proposition 2.2. By Corollary
2.1, the position $x_0(y_a)$ of the singularity (near $b$) of the
solution $y(y_a,x)$, analytically continued along $s$, is a
regular function of $y_a$. By Proposition 2.2, $u_b={1\over y_b}$,
where $y_b=y(y_a,b)$, as a function of $x_0$, is given by \be
u=u_b =\sum_{k=1}^{\infty} \beta_k (x_0-b)^{k\over 2}, \ \ \beta_1
\neq 0.\ee Substituting into this expression a regular function
$x_0(y_a), \ x_0(y^0_a)=b,$ we get \be u=u_b =\sum_{k=1}^{\infty}
c_k (y_0-y^0_a)^{k\over 2}, \ \ c_1 \neq 0.\ee Finally, expressing
$y_b={1\over u_b}$ through $u_b$ via (4.3), we get the required
formula (4.1).

\section{Analysis of a local model}
\setcounter{equation}{0}
Here, we want to investigate the local behaviour of the solutions
near a generic fixed singularity. That is to say, we assume that
the polynomial $q(x)$ has a simple zero and we indeed replace the
equation by the following:
\be \frac{dy}{dx}=cxy^{3}+y^{2}.\ee The
change of unknown function
$$y=\frac{v}{x},$$
yields
$$-\frac{v}{x^{2}}+\frac{v'}{x}=
\frac{cv^{3}}{x^{2}}+\frac{v^{2}}{x^{2}},$$  and thus: \be
v'=\frac{1}{x}[cv^{3}+v^{2}+v], \ee which obviously separates.
\vskip 1pt Write:
$$\frac{1}{cv^{3}+v^{2}+v}=\frac{\alpha}{v}+\frac{\beta}{v-v_{1}}
+\frac{\gamma}{v-v_{2}},$$ with:
$$v_{1}=\frac{-1+\sqrt{1-4c}}{2c},
\quad v_{2}=\frac{-1-\sqrt{1-4c}}{2c},$$ and
$$\alpha=1, \quad \beta=\frac{1}{cv_{1}(v_{1}-v_{2})},
\quad \gamma=\frac{1}{cv_{2}(v_{2}-v_{1})}.$$ Integrating equation
(5.2) we get for each its solution $v(x)$
$$v(v-v_{1})^{\beta}(v-v_{2})^{\gamma}=K\cdot x,$$ for a certain
constant $K$. Equivalently
$$y(xy-v_{1})^{\beta}(xy-v_{2})^{\gamma}=K,$$
or \be
y(1-\frac{xy}{v_{1}})^{\beta}(1-\frac{xy}{v_{2}})^{\gamma}={K\over
{v_1^{\beta}v_2^{\gamma}}}=K'.\ee Notice that the only ``fixed
singularity" of the equation (5.1) is $x=0$. To start with, let us
take this point $x=0$ as the initial point $a$. Now, the constant
$K'$ in (5.3) is evaluated by setting $x=0$ and $y=y_{0}$ and this
yields $K'=y_{0}$. Therefore, the solution $y(y_0,x)$ of (5.1)
satisfying $y(y_0,0)=y_0$ is given by \be
y(x)(1-\frac{xy(x)}{v_{1}})^{\beta}
(1-\frac{xy(x)}{v_{2}})^{\gamma}=y_0.\ee Substituting into (5.4)
the point $x=b$, we get the relation between $y_b=y(b)$ and
$y_{0}$ in the form: \be
y_{0}=y(1-\frac{by_b}{v_{1}})^{\beta}(1-\frac{by_b}{v_{2}})^{\gamma}.
\ee Now, we are interested in the ``limit cycles" of the equation
(5.1), i.e. in its solutions $y(x)$ satisfying $y(0)=y(b)$. This
relation together with (5.5) gives the equation for the limit
cycles, which are (besides the solution $y\equiv 0$ of (5.1)) in
one-to-one correspondence with the solutions of: \be
(1-\frac{by_{0}}{v_{1}})^{\beta}(1-\frac{by_{0}}{v_{2}})^{\gamma}=1.
\ee At this point we have to clarify the geometric interpretation
of the ``limit cycles", as defined above. The problem is that the
solutions of (5.1) are multivalued functions. The accurate
interpretation of the equation (5.6) is that {\it the algebraic
curve $Y=Y_{y_0}$, defined by (5.4), passes through the points
$(0,y_0)$ and $(b,y_0)$.} Certainly, this curve $Y$, parametrized
as $y=y(x)$, satisfies differential equation (5.1). But a priori
we do not even know whether $Y$ is connected. So (5.6) by itself
does not exclude a possibility that the points $(0,y_0)$ and
$(b,y_0)$ belong to different leaves of the solutions of the
differential equation (5.1).

\medskip

Below we show that in fact for $y_0\neq 0$ the curve $Y=Y_{y_0}$
is connected. This allows us to give the following interpretation
to the equation (5.6): {\it for each $y_0$ satisfying (5.6) there
exists a path $s$ from $0$ to $b$ such that the solution
$y(y_0,x)$ can be analytically continued along $s$, and this
continuation satisfies $y(y_0,b)=y_0$.}

\medskip

We now choose some specific values for the free parameter $c$ in
order to bring some light on possible solutions of (5.6). Assume
that
$$1-4c=\delta^{2},$$
where $\delta=-{1\over {2n+1}}$ with $n$ an integer. In that case
we obtain:
$$\beta=n, \ \gamma= -n-1,$$ (see Section 6 below for more detailed
computations). Limit cycles of equation (5.1), in the
interpretation given above, are in correspondence with the
solutions of: \be
(1-\frac{by_{0}}{v_{1}})^n=(1-\frac{by_{0}}{v_{2}})^{n+1}. \ee One
can easily show that for large integer values of $n$ equation
(5.7) has $n$ distinct complex solutions $y_0^j, \ j=1,\dots, n$.
(see Section 7 below). Consider the local solutions
$y^j(x)=y(y_0^j,x)$ at the origin satisfying $y^j(0)=y_0^j.$
Combining equation (5.7) with Theorem 6.1 below which describes
the monodromy of the solutions of (5.1) we get the following
result: \bt For $c={1\over 4}(1- {1\over {(2n+1)^2}})$ and $b\neq
0$ equation (5.1) has $n$ different ``limit cycles", i.e. local
solutions $y^j(x)$ at the origin, $j=1,\dots,n,$  and paths $s^j$
from $0$ to $b$, such that each $y^j(x)$ being analytically
continued along $s^j$ satisfies $y(0)=y(b)$.\et The proof of this
theorem is given at the end of Section 6 below. From the
description given in Section 6 it follows that the paths $s^j$
have the following form: $s^j$ goes from zero to the (only)
singularity of $y^j(x)$, turns once around this singularity,
returns to zero, makes $m\leq n$ turns around zero, and finally
comes to $b$.

\medskip

Therefore, in this example we see that the equation (5.1) may have
as many complex limit cycles as we wish, when $n$ is increased,
although the degree of the coefficients of this equation remains
bounded. This phenomenon reminds (in much simpler setting) the
counterexample due to Yu. Iliashenko of Petrowski-Landis claim
(\cite{Ily}). It should also be compared with the examples of
differentiable Abel equations discussed by A. Lins Neto
(\cite{Ln}). Note that Khovanski fewnomials theory, or, rather,
``additive complexity" arguments (see \cite{Kho,Ris}) imply that
the number of real roots of equation (5.7) remains bounded
independently of $n$. So that this example does not provide a
counterexample to real Hilbert-Pugh problem.

\medskip

\noindent{\bf Remark.} One can investigate the situation for
another choice of the parameter $c$. If we put ${n\over {n+1}}=k$
and let $k$ to be a large integer, this corresponds to $n\approx
-1$ or $c\approx 0$. The equation (5.4) takes the form \be
(1-\frac{xy}{v_{1}})^k=y_0^{\mu}(1-\frac{xy}{v_{2}}). \ee while
the equation (5.7) takes the form \be
(1-\frac{by_{0}}{v_{1}})^k=(1-\frac{by_{0}}{v_{2}}). \ee The
investigation of this case may be important since as $c$ tends to
zero, equation (5.1) tends to the integrable equation $y'=y^2$.
\medskip
We now consider the case $c=\frac{1}{4}$ and the equation
$$\frac{dy}{dx}=\frac{1}{4}xy^{3}+y^{2}.$$
With $y=v/x$, this equation yields:
$$\frac{dv}{dx}=\frac{1}{4x}(v^{3}+4v^{2}+4),$$
which separates and gives the solution $y(x)$ corresponding to the
initial data $y_{0}$ as the solution to:
$$\frac{y}{xy+2}{\rm e}^{\frac{2}{xy+2}}=\frac{{\rm e}}{2}y_{0}.$$
Periodic orbits correspond to solutions of $y(1)=y(0)=y$ to
$$\frac{2}{y+2}{\rm e}^{\frac{2}{y+2}}={\rm e}.$$
If we change
$y=\frac{2\xi}{1-\xi},$
this yields
$$1-\xi={\rm e}^{\xi}.$$
We write $\xi=x+{\rm i}y,$
and derive the two equations
$$1-x={\rm e}^{x}{\rm cos}y,$$
$$-y={\rm e}^{x}{\rm sin}y.$$
Note that if $(x,y)$ is a solution, then $(x,-y)$ is also a
solution. Then we can assume $y>0$. Second equation implies ${\rm sin}(y)<0$
and we restrict ourselves now to $y\in ](2n+1)\pi, (2n+2)\pi[$.
Now we plug $x=-{\rm log}(-\frac{{\rm sin}y}{y})$ into the first
equation. This displays:
$$F(y)=1+{\rm log}(-\frac{{\rm sin}y}{y})+\frac{y}{{\rm tan}y}=0.$$
Then we note that as $y\to (2n+1)\pi$, $F(y)\to +\infty$ and that
as $y\to (2n+2)\pi$, $F(y)\to -\infty$. There is thus at least one
solution (and in fact a single one) in the interval. The Abel
equation has thus infinitely many limit cycles.

\medskip

The basic example (5.1) can be used to generate a family of
similar examples by composition. Composition appears quite
naturally in the subject (see both the classics (Abel,Liouville,...) and more
recent contributions (\cite{Bfy1}-\cite{Bfy5}). Consider the Abel
equations of the form \be y'=cP(x)p(x)y^{3}+p(x)y^{3}, \ee where
$p(x)$ is an arbitrary polynomial, and $P(x)$ is the
anti-derivative of $p(x)$ which vanishes at $x=0$. The change of
variables $w=P(x)$ brings (5.9) to the form \be {dy\over
{dw}}=cwy^3+y^2.\ee Applying the above given analysis of this last
equation, we see that the solution $y(x)$ to the equation (5.9)
satisfying $y(0)=y_{0}$, solves the implicit algebraic equation:
\be y(1-\frac{P(x)y}{v_{1}})^{\beta}
(1-\frac{P(x)y}{v_{2}})^{\gamma}=y_{0}. \ee Hence also the limit
cycles of (5.9) can be investigated in a similar way. Notice,
however, that a special composition structure of the solutions of
(5.9), namely, that each its solution $y$ can be represented as
$y(x)=\tilde y(P(x))$, for $\tilde y(w)$ solving (5.10), implies
the following: for each $a,b$ with $P(a)=P(b)$ we have $y(a)\equiv
y(b)$.

\section{Ramification of solutions of $\frac{dy}{dx}=cxy^{3}+y^{2}$}
\setcounter{equation}{0}

To study in detail the ramification of the solutions of the Abel
equation (5.1), $\frac{dy}{dx}=cxy^{3}+y^{2}$, we choose the
parameter $c$ in this equation in the same way as above. We would
like $c$ to tend to $1\over 4$, which is a ``discriminant point"
for the denominator $cv^{3}+v^{2}+v$ appearing after separation of
variables in (5.1). On the other hand, we want the first integral
to remain algebraic. So let us write
$$1-4c=\delta^{2}, \ \ c={1\over 4} - {\delta^{2}\over 4}$$ where
$\delta$ is small. In this case we obtain:
$$v_{1}=-\frac{2}{1+\delta}\approx -2+2\delta,
\quad v_{2}=-\frac{2}{1-\delta}\approx -2-2\delta.$$ For $\beta$
and $\gamma$ we get, respectively, $$ \beta=-{{1+\delta}\over
{2\delta}}, \ \ \gamma={{1-\delta}\over {2\delta}}.$$ Notice, that
$\beta+\gamma=-1.$ Thus, taking $\delta=-{1\over {2n+1}}$ for $n$
an integer, we obtain $\beta=n, \ \ \gamma=-n-1.$ Let us fix this
choice of $c$.

\medskip

The first integral (5.4) takes now the form \be
y_0={\frac{y(1-{{xy}\over v_1})^n}{(1-{{xy}\over
v_2})^{n+1}}}=H(x,y).\ee In the rest of this section we
investigate in some detail the solutions of the algebraic equation
(6.1). First of all, we notice that for $y_0=0$ we indeed get {\it
two separate leaves}: the straight line $Y_0^0 =\{y\equiv 0\}$ and
the hyperbola $Y_0^1 =\{y={v_1\over x}\}$. As the zero curves of
$H(x,y)$ the first has the multiplicity $1$, while the second has
the multiplicity $n$.

The second remark is that the pole locus of $H(x,y)$ is the
hyperbola $Y_{\infty} =\{y={v_2\over x}\}$ which has the
multiplicity $n+1$. As we can expect, both the hyperbolas above
are solutions of the differential equation (5.1). This can be
checked by a direct substitution.

\subsection{Critical points of H(x,y)}

Let us find the critical points of the function $H(x,y)$. \bl All
the critical points of the function $H(x,y)$ (in fact, all the
points with $H_y(x,y)=0$) are situated on the hyperbola $Y_0^1
=\{y={v_1\over x}\}$.\el \pr After differentiating $H$ with
respect to $x$ and $y$, cancelling the common degrees of
$(1-{{xy}\over v_2})$, equating the numerator to zero, and some
computations, using, in particular, the identities
$${n\over v_2}-{{n+1}\over v_1}=0, \ \
{n\over v_1}-{{n+1}\over v_2}=1,$$ we obtain the following system
of equations: $$(1-{{xy}\over v_2})^{n}H_y=(1-{{xy}\over
v_1})^{n-1}=0,$$ $$(1-{{xy}\over v_2})^{n}H_x=(1-{{xy}\over
v_1})^{n-1}(1-{{xy}\over {v_1v_2}})=0.$$ The common zeroes of this
system lie exactly on the parabola $Y_0^1 =\{y={v_1\over x}\}$.
Notice that the partial derivative $H_x$ vanishes, in addition, on
the hyperbola $y={{v_1v_2}\over x}.$ Hence, the points of this
hyperbola are zeroes of the derivative $y'$ of the solutions of
(5.1) passing through these points. Of course, this can be checked
by the direct substitution.

\medskip

\noindent{\bf Remark 1.} The fact that $H$ is the first integral
of the equation (5.1), and hence its level curves must be locally
graphs of a regular function at each finite point, does not
exclude by itself possible critical points of $H$ - compare the
points of the hyperbola $Y_0^1$.

\medskip

\noindent{\bf Remark 2.} Instead of the rational equation (6.1) we
can consider the equivalent polynomial equation \be y(1-{{xy}\over
v_1})^n -y_0(1-{{xy}\over v_2})^{n+1}=0.\ee The advantage of (6.1)
is that the initial value $y_0$ appears there just as the right
hand side.

Now, differentiating (6.2) we get the following system of
equations:

$$(1-{{xy}\over v_1})^n + {{nxy}\over v_1}(1-{{xy}\over
v_1})^{n-1} - {{(n+1)xy_0}\over v_2}(1-{{xy}\over v_2})^{n}=0,$$
$${{ny^2}\over v_1}(1-{{xy}\over
v_1})^{n-1} - {{(n+1)yy_0}\over v_2}(1-{{xy}\over v_2})^{n}=0.$$
Multiplying the first equation by $y$ and the second by $x$ and
taking the difference, we get $$y(1-{{xy}\over v_1})^n=0.$$ So
either $y=0$ or the point $(x,y)$ belongs to $Y_0^1$. If $y=0$ the
second equation above is satisfied, while the first equation gives
$x={v_2\over {(n+1)y_0}}.$ So the equation (6.2) has an additional
critical point, not on the hyperbola $Y_0^1$. Notice, however,
that for any $y_0\neq 0$ this point does not belong to the
solution curve of (6.2), while for $y_0=0$ it is at infinity.

\subsection{Singularities of solutions of
$\frac{dy}{dx}=cxy^{3}+y^{2}$}

The only ``fixed" singularity of the equation (5.1) is the origin
$x=0$. Let us start with the ``movable" singularities $x_0\neq 0$
of the solutions (compare with the general results of Section 2).
\bp For $y_0\neq 0$ the solution $y(y_0,x)$ has the only movable
singularity at the point $x_0(y_0)={\kappa \over y_0}$, where
$\kappa = -v_2({v_2\over v_1})^n=-{2\over
{1-\delta}}({{1-\delta}\over {1+\delta}})^n \approx -2e$ for
$\delta$ small. Exactly one local branch of $y(y_0,x)$ takes an
infinite value and has a ramification of order $2$ at $x_0(y_0)$,
while the other $n-1$ local branches are regular at this point and
take there $n-1$ different finite values.\ep \pr Denoting, as
above, $1\over y$ by $u$ we get from (6.1) \be y_0={\frac
{(u-{x\over v_1})^n}{(u-{x\over v_2})^{n+1}}}. \ee Substituting
here $u=0$ (and assuming $x\neq 0$ and so cancellation is
possible) we get $x_0(y_0)={\kappa \over y_0}$. To get the series
expansions we rewrite (6.3) as follows: \be x={\kappa \over
y_0}{\frac {(1-{uv_1\over x})^n}{(1-{uv_2\over x})^{n+1}}}={\kappa
\over y_0}(1+{u\over x}[(n+1)v_2-nv_1] + A({u\over
x})^2+\dots).\ee Since $(n+1)v_2-nv_1=0$, we can rewrite the last
expression as \be x={\kappa \over y_0}(1+A({u\over x})^2+\dots), \
\ x-{\kappa \over y_0}=B({u\over x})^2(1+\dots).\ee This shows
that $x-{\kappa \over y_0}$ as a function of $u$ has a second
order zero at $u=0$, and hence $u$ as a function of $x$ has at
$x_0={\kappa \over y_0}$ a second order branching. (We do not
prove that the coefficient $B$ above is different from zero, since
this fact was shown in general form in Section 2 above). This
completes the description of the branch passing through the point
$(x_0(y_0),\infty)$.

\medskip

Each other branch of the solution $y(y_0,x)$ (i.e. of the curve
$H(x,y)=y_0$) over $x_0$ is regular, since by Lemma 6.1 all the
singularities of $H$ belong to the level curve $H=0$. Since for
any fixed $x$ the total number of the solutions of $H(x,y)=y_0$
with respect to $y$, counted with multiplicities, is $n+1$, and
since the multiplicity of the singular branch is $2$, there are
exactly $n-1$ regular local branches of the curve $H(x,y)=y_0$
over $x_0$. This completes the proof of the proposition.

\medskip

\noindent{\bf Remark.} Exactly as in Sections 3 and 4 above, we
can use the series (6.5) to analyze the local structure of
singularities of the Poincar\'e mapping. Indeed, for a fixed
$x=x_0$ we can rewrite (6.5) as \be y_0-{\kappa \over
x_0}=D({u\over x_0})^2(1+\dots), \ee and we get a second order
zero of $y_0-{\kappa \over x_0}$ as a function of $u$ and a second
order ramification of $u$ as a function of $y_0$.

\medskip

The next step is to investigate the structure of the fixed
singularity $x=0$. \bp For $y_0\neq 0$ the solution $y(y_0,x)$ has
over $x=0$ two local components: the regular one, passing through
the point $(0,y_0)$, and the singular one, passing through the
point $(0,\infty)$. The singular component is represented by the
Puiseux series \be{1\over {y(x)}}= u(x)={1\over v_1}x+\mu
x^{1+{1\over n}}+\dots,\ee with $$\mu = y_0^{1\over n}({1\over
{2n+1}})^{1+{1\over n}} \approx {1\over {2n}}y_0^{1\over n}.$$ In
particular, the local monodromy acts as a cyclic permutation of
the infinite branches.\ep \pr Let us rewrite the equation (6.3) in
the form \be y_0(u-{x\over v_2})^{n+1}=(u-{x\over v_1})^n. \ee We
have to find the Puiseux expansion of the curve given by (6.8) at
the point $(u,x)=(0,0)$. To simplify the presentation, we use the
following ``einsatz": \be u(x)={1\over v_1}x+\mu x^{\nu}
+\dots.\ee Substituting (6.9) to (6.8) we get $$ y_0[({1\over
v_1}-{1\over v_2})x+\mu x^{\nu}+\dots]^{n+1}=\mu^n x^{\nu
n}+\dots.$$ Comparing the leading degrees and coefficients, we
obtain $$ y_0({1\over v_1}-{1\over v_2})^{n+1}x^{n+1}+\dots=\mu^n
x^{\nu n}+\dots,$$ and hence $$\mu=y_0^{1\over n}({1\over
v_1}-{1\over v_2})^{{n+1}\over n}=y_0^{1\over
n}(-\delta)^{1+{1\over n}}\approx {1\over {2n}}y_0^{1\over n}, \ \
\nu={1+{1\over n}}.$$

\subsection{Global ramification of solutions}

According to Proposition 6.1, there are only two singularities of
the solution $y(y_0,x)$: the fixed singularity at $x=0$ and the
movable singularity at $x_0(y_0)={\kappa \over y_0}$. The original
local branch at $x=0$ of $y(y_0,x)$ is regular at the origin.
Hence, it can be analytically extended as a regular univalued
function into the disk ${\cal D}=D_{\vert{\kappa \over
y_0}\vert}$, centered at $x=0$. \bl The regular branch of
$y(y_0,x)$ on the disk $\cal D$ has a singularity at the boundary
point ${\kappa \over y_0}$.\el \pr Take $y_0$ positive. By our
choice of the parameter $c\approx {1\over 4}$ we have $c>0$.
Therefore, the right hand side of (5.1) is positive, and bounded
from below by $y^2$, and hence its solution blows up in finite
time on the semi-axis $x>0$. By Proposition 6.1, this happens
exactly at the point $x_0(y_0)={\kappa \over y_0}$. This proves
Proposition 6.2 for $y_0$ positive. Now, as $y_0$ moves along the
circle $\vert y_0 \vert=Const$, the singularity $x_0(y_0)={\kappa
\over y_0}$ of the regular univalued function $y(y_0,x)$ on the
disk $\cal D$ moves along the boundary of this disk. Since
$y(y_0,x)$ analytically depends on $y_0$, the point $x_0(y_0)$
remains its singularity. This completes the proof.

\medskip

Let the value $y_0\neq 0$ be fixed. Consider the loop $\omega$
following the straight segment from $0$ to the singular point
$x_0(y_0)$, then going around this point in a counter-clockwise
direction along a small circle, and then returning to $0$ along
the same straight segment. \bl The regular branch at $x=0$ of the
solution $y(y_0,x)$ analytically continued along the loop
$\omega$, returns at $x=0$ to one of the infinite branches of the
solution.\el \pr Since $y(y_0,x)$ has a second order ramification
at $x_0$, after one turn around this point we get {\it another}
branch of the solution. As we return to zero, we stay on this new
branch, different from the initial (regular) one. But by
Proposition 6.2 all the branches, except the initial one, tend to
$\infty$ at $x=0$.

\medskip

Now we have enough tools to prove one of the main properties of
the solutions of (5.1), as given by the first integral
$H(x,y)=y_0$: \bt For each $y_0\neq 0$ the solution curve
$Y_{y_0}=\{H(x,y)=y_0 \}$ is irreducible. The analytic
continuation of the local solution $y(y_0,x)$ at zero along the
loop $\omega$ and then several turns around zero transform this
local branch to each one of the $n$ remaining branches of
$Y_{y_0}$.\et \pr By Lemma 6.3 continuation along $\omega$
transforms the local regular branch of $y(y_0,x)$ at zero into one
of the infinite branches. By proposition 6.2, each turn around
zero results in a cyclic permutation of the $n$ infinite branches.
Hence, in at most $n$ turns each other infinite branch can be
obtained.

\medskip

\noindent{\bf Remark.} Another proof of Theorem 6.1 can be
obtained by computing the ramification of $Y_{y_0}$ at $x=\infty$.
Rewriting the equation (6.3) in the form $$ y_0 ({1\over
y}-{x\over v_2})^{n+1}=({1\over y}-{x\over v_1})^n, $$ and then
substituting $x={1\over w}$, we obtain \be \kappa
wy(y-v_1w)^n=y_0(y-v_2w)^{n+1}. \ee The einsatz $y=v_2w+\eta
w^{\rho}+\dots$ leads to $$\kappa w(v_2w+\eta
w^{\rho}+\dots)[(v_2-v_1)w+\eta w^{\rho}+\dots]^{n}=y_0(\eta
w^{\rho}+\dots)^{n+1},$$ which produces, via comparing the leading
terms, $$\kappa v_2(v_2-v_1)^n w^{n+2}=y_0 \eta^{n+1} w^{\rho
(n+1)} $$ and $$ \eta=({{\kappa v_2} \over y_0})^{1\over
{n+1}}(v_2-v_1)^{n\over {n+1}}, \ \ \rho=1+{1\over {n+1}}.$$ We
conclude that all the $n+1$ branches of $Y_{y_0}$ tend to zero at
$w=0$ or $x=\infty$, and that the local monodromy around infinity
produces a cyclic permutation of these $n+1$ branches.

\medskip

\noindent{\bf Proof of Theorem 5.1.} Let us remind equation (5.7):
$$ (1-\frac{by_{0}}{v_{1}})^n=(1-\frac{by_{0}}{v_{2}})^{n+1}. $$
One has to show that for large integer values of $n$ equation
(5.7) has $n$ distinct complex solutions $y_0^j, \ j=1,\dots, n$.
Consider the local solutions $y^j(x)=y(y_0^j,x)$ at the
origin satisfying $y^j(0)=y_0^j.$ The analytic continuation of
$y^j(x)$ gives an algebraic curve $Y_{y_0^j}$ satisfying equation
(5.4) with $y_0=y_0^j$. Now equation (5.7) says exactly that the
points $(0,y_0^j)$ and $(b,y_0^j)$ belong to $Y_{y_0^j}$. By
Theorem 6.1 this curve is irreducible, and we can pass from the
point $(0,y_0^j)$ to the point $(b,y_0^j)$ via the analytic
continuation of the local branch $y^j(x)=y(y_0^j,x)$ at the origin
along the path $s^j$, as described in Theorem 5.1. This completes
the proof.

\bibliographystyle{amsplain}

\end{document}